\documentclass[11pt]{article}
\usepackage{amssymb}
\usepackage{amsmath}
\usepackage{theorem}
\usepackage{latexsym}
\reversemarginpar
\theoremheaderfont{\bf} 
\theorembodyfont{\sl}
\newtheorem{theo}{Theorem}[section]
{\theorembodyfont{\rm} \newtheorem{defi}[theo]{Definition}}
{\theorembodyfont{\rm} \newtheorem{exa}[theo]{Example}}
{\theorembodyfont{\rm} \newtheorem{rem}[theo]{Remark}}
\newtheorem{prop}[theo]{Proposition}
\newtheorem{cor}[theo]{Corollary}
\newtheorem{lemma}[theo]{Lemma}
{\theorembodyfont{\rm}}
{\theorembodyfont{\rm}}
\newenvironment{proof}{{\sc Proof:}}{\mbox{}\hfill$\Box$\par}
\newcommand{\eqnref}[1]{~\mbox{$(${\rm \ref{#1}}$)$}}

\newcommand{\comment}[1]{\marginpar[{#1}]{{#1}}}
\newcommand{\junk}[1]{}
\newcommand{\DS}{\displaystyle}

\newcommand{\N}{{\mathbb N}}
\newcommand{\F}{{\mathbb F}}

\newcommand{\cC}{{\mathcal C}}

\newcommand{\cI}{{\mathcal I}}
\newcommand{\cS}{{\mathcal S}}

\newcommand{\CC}{convolutional code}
\newcommand{\CCC}{cyclic convolutional code}
\newcommand{\ve}[1]{\mbox{$\varepsilon^{(#1)}$}}
\newcommand{\scy}{\mbox{$\sigma$-cyclic}}
\newcommand{\rank}{\mbox{\rm rank}\,}
\newcommand{\spann}{\mbox{\rm span}\,}
\newcommand{\equivs}{\equiv_{_{\!\sigma}}}
\newcommand{\nequivs}{\not\equiv_{_{\!\sigma}}}
\newcommand{\AutF}{\mbox{${\rm Aut}_{\mathbb F}$}}
\newcommand{\im}{\mbox{\rm im}\,}
\renewcommand{\mod}{\mbox{\rm mod}\,}

\newcommand{\nkd}{\mbox{$(n,k,\delta)$}}
\newcommand{\dist}{\mbox{\rm dist}}
\newcommand{\wt}{\mbox{\rm wt}}
\newcommand{\id}{\mbox{\rm id}}

\newcommand{\ideal}[1]{\mbox{$\langle{#1}\rangle$}}
\newcommand{\Flaurent}{\mbox{$\F(\!(z)\!)$}}
\newcommand{\Azs}{\mbox{$A[z;\sigma]$}}

\newcommand{\p}{\mbox{$\mathfrak{p}$}}
\renewcommand{\v}{\mbox{$\mathfrak{v}$}}

\newcommand{\lideal}[1]{\mbox{$^{^{\bullet\!\!}}\langle{\, #1\, }\rangle$}}

\newcommand{\db}[1]{\mbox{$[\![#1]\!]$}}

\newcounter{abc}
\newcounter{def}
\newenvironment{romanlist}{\begin{list}{(\roman{abc})\hfill}{\usecounter{abc}
     \topsep0ex \labelwidth.7cm \leftmargin.7cm \labelsep0cm
     \rightmargin0cm \parsep0ex \itemsep.6ex
     \partopsep0ex}}{\end{list}}
\newenvironment{alphalist}{\begin{list}{(\alph{abc})\hfill}{\usecounter{abc}
     \topsep0ex \labelwidth.7cm \leftmargin.7cm \labelsep0cm
     \rightmargin0cm \parsep0ex \itemsep.6ex
     \partopsep0ex}}{\end{list}}
\newenvironment{arabiclist}{\begin{list}{(\arabic{abc})\hfill}{\usecounter{abc}
     \topsep0ex \labelwidth.7cm \leftmargin.7cm \labelsep0cm
     \rightmargin0cm \parsep0ex \itemsep.6ex
     \partopsep0ex}}{\end{list}}

\title{On the Parameters of Convolutional Codes with Cyclic Structure}
\date\today
\author{Heide Gluesing-Luerssen\thanks{
           Department of Mathematics, University of Kentucky, 715
           Patterson Office Tower, Lexington, KY 40506, USA; heidegl@ms.uky.edu}
        \mbox{}\
        and Barbara Langfeld\thanks{
        Kombinatorische Geometrie (M9),
        Zentrum Mathematik, Technische Universit\"at M\"unchen, Boltzmannstr.~3,
        85747 Garching bei M\"unchen, Germany;
        langfeld@ma.tum.de}
        }
\begin{document}
\maketitle

\begin{abstract}
\noindent
In this paper convolutional codes with cyclic structure
will be investigated.
These codes can be understood as left principal ideals in a suitable
skew-polynomial ring.
It has been shown in~\cite{GS02} that only certain combinations of the parameters
(field size, length, dimension, and Forney indices) can occur for cyclic codes.
We will investigate whether all these combinations can indeed be realized by a suitable
cyclic code and, if so, how to construct such a code.
A complete characterization and construction will be given for minimal cyclic codes.
It is derived from a detailed investigation of the units in the
skew-polynomial ring.
\end{abstract}



\section{Introduction}
\setcounter{equation}{0}

The two most important classes of codes used in practice are block codes and
convolutional codes.
While both classes play an equally important role in engineering practice, the
theory of \CC{}s is much younger and not nearly as developed as the theory of
block codes.
The foundation of the mathematical theory of \CC{}s was laid only in the seventies of the last
century by the articles of Forney, see e.~g.~\cite{Fo70}.
It led to quite some mathematical investigation in that decade among which
are basically two groups of papers.

The first group~\cite{MCJ73,Ju73,Ju75} deals with the construction of
convolutional codes with large distance, mainly by using cyclic block codes and
resorting to the weight-retaining property for bridging the gap between cosets
of polynomials in the block code case and vector polynomials in the convolutional case.
These ideas were resumed later again in~\cite{SGR01}, leading to the
construction of MDS convolutional codes.

The second group of papers~\cite{Pi75,Pi76,Ro79} initiated a completely different approach.
In the paper~\cite{Pi76} it was investigated for the first time as to how
cyclic structure has to be understood for a convolutional code itself.
The first crucial fact being found was that cyclic structure in the classical
sense (i.~e.\ invariance under the cyclic shift) is not an appropriate concept for convolutional codes.
Precisely, it was shown in~\cite{Pi76} that each convolutional code, that is invariant
under the cyclic shift, has complexity zero, hence is a block code.
This insight has led Piret to a different, much more complex notion of cyclicity, which then
was further generalized by Roos~\cite{Ro79}.
In the simplest form, this structure can be understood as a sort of graded shift in the
coefficients of the polynomial codewords.
The precise notion will be given in Section~2.
At this point we only want to mention that cyclic codes of length~$n$ over the field~$\F$ can be
understood as certain left ideals in a skew-polynomial ring $A[z;\sigma]$, where
$A=\F[x]/_{\ideal{x^n-1}}$, the variable~$z$ represents the delay operator, and~$\sigma$
determines the non-commutative structure.
Both Piret and Roos gave several examples of convolutional codes, that are cyclic in this
new sense.
They also computed (or estimated) the distances which turned out to be very good.

Although these papers initiated an algebraic theory of cyclic convolutional codes,
they did not come very far and the topic came to a halt.
Only recently it has been resumed in~\cite{GS02}.
Therein an algebraic theory of \CCC{}s, fully in terms of ideals in the skew-polynomial ring, has been
established.
It leads to a nice, yet nontrivial, generalization of the algebraic theory of cyclic
block codes.
The translation from ideals into polynomial vectors is achieved by suitable circulant matrices.
In particular, \CCC{}s are principal left ideals (thus have a generator polynomial),
they are also left annihilators of right ideals (thus have a parity check polynomial),
the parameters can be computed in terms of these polynomials,
and the dual of a cyclic code is cyclic again.
Moreover, in~\cite{GS03} plenty of examples of \CCC{}s are given, their distances are all
optimal in the sense that they attain the Griesmer bound.
All this indicates that the notion of cyclicity as introduced by Piret is the
appropriate one for \CC{}s not only when it comes to the algebraic theory, but also for
constructing good codes.

In this paper we will continue the algebraic theory as it was set up in~\cite{GS02}.
It is a consequence of the results in~\cite{GS02} that only certain combinations of parameters
(field size, length, dimension, and Forney indices) can occur for cyclic codes; see also
Theorem~\ref{T-ideals}(4) below.
We seek to investigate whether all these combinations do really occur.
The key role for this aim is played by so called minimal \CCC{}s, these are cyclic codes that have no
proper cyclic subcodes.
They form the building blocks of all cyclic codes in the sense that each cyclic code is
the direct sum of minimal codes and the Forney indices of the code are given by the union
of the Forney indices of each component.
Minimal codes have a very simple ideal theoretic description in terms of their generator polynomial,
see Proposition~\ref{P-minC}.
Moreover, for these codes all Forney indices are the same, hence these codes are compact in
the sense of~\cite[Cor.~4.3]{McE98}.
This makes these codes also very important from a coding point of view since
compact codes are in general good candidates for having a large distance.
(for instance codes attaining the generalized Singleton bound are always compact,
see~\cite{RoSm99}).
We will show that under a certain necessary and sufficient condition any arbitrarily chosen
Forney index can be realized by a suitable minimal cyclic code and we will show how to construct
such a code.
This result will then be further exploited for investigating non-minimal codes with prescribed
Forney indices.

The outline of the paper is as follows.
The end of the introduction is devoted to the basic notions of convolutional coding theory.
Thereafter in Section~2 we will introduce cyclicity for convolutional codes
along with the algebraic machinery and the main results from~\cite{GS02} as needed for our
purposes.
In Section~3 we turn to minimal \CCC{}s.
Their investigation amounts basically to a detailed study of the units in the skew
polynomial ring $\Azs$.
This will lead us to the existence of minimal codes with prescribed Forney indices under certain
necessary and sufficient conditions.
Finally, in Section~4 we will turn to certain direct sums of minimal codes.
These direct sums are specific in the sense that the generator polynomials of the minimal
components are pairwise orthogonal, resulting in an easy handling of the direct sum.
The existence result from Section~3 will be extended to these codes.

We will end the introduction with the basic notions of convolutional coding theory.
Convolutional codes are certain submodules of~$\F[z]^n$, where~$\F$ is a finite
field.
Before presenting the definition we wish to recall that each submodule~$\cS$ of
$\F[z]^n$ is free and therefore can be written as
\[
   \cS=\im G:=\big\{uG\,\big|\, u\in\F[z]^k\big\}
\]
where~$k$ is the rank of~$\cS$ and $G\in\F[z]^{k\times n}$ is a matrix
containing a basis of~$\cS$. Any such matrix~$G$ is called a {\em
generator matrix\/} of the module~$\cS$. It is unique up to left
multiplication by a unimodular matrix, that is, for any pair of
matrices $G,\,G'\in\F[z]^{k\times n}$ having full row rank the
identity $\im G=\im G'$ is equivalent to $G'=VG$ for some matrix
$V\in Gl_k(\F[z])$.
This makes the following notions well-defined.

\begin{defi}\label{D-CC}
Let $\F$ be any finite field and let
$G\in\F[z]^{k\times n}$ be a matrix of rank~$k$.
\begin{alphalist}
\item The number
      $\delta:=\delta(G):=\max\{\deg\gamma\mid \gamma$ is a
      $k$-minor of $G\}$ is called the {\em complexity\/} of the submodule
      $\im G$ or of the matrix~$G$.
\item The submodule $\cC:=\im G\subseteq\F[z]^n$ is called a  {\em
      convolutional code over $\F$ with parameters\/} $(n,k,\delta)$ if
      it has complexity~$\delta$ and
      the matrix $G$ is right invertible, i.~e. if
      there exists some matrix $\tilde{G}\in\F[z]^{n\times k}$ such that
      $G\tilde{G}=I_k$.
      In this case the parameter~$n$ is called the {\em length\/} of the code.
\end{alphalist}
\end{defi}
Since every right invertible matrix $G\in\F[z]^{k\times n}$ can be
completed to a unimodular matrix (e.g.\ by using the Smith normal form)
one has the following properties.

\begin{rem}\label{R-dirSummand}
\begin{alphalist}
\item The convolutional codes over~$\F$ of length~$n$ are the direct summands
      of the module $\F[z]^n$.
\item Each convolutional code $\cC\subseteq\F[z]^n$ has a parity check matrix, that
      is, there exists a matrix $H\in\F[z]^{n\times(n-\text{rk}\,\cC)}$ such that
      $\cC=\ker H:=\{v\in\F[z]^n\mid vH=0\}$.
\end{alphalist}
\end{rem}
Part~(b) can be considered as the main reason for restricting to direct summands
rather than arbitrary submodules for \CC{}s.
A parity check matrix is an important tool for data transmission, it is needed for
checking whether or not the received data are erroneous.

The following property of \CC{}s will be needed later on.
\begin{lemma}\label{L-directsum}
Let $\cC,\,\hat{\cC}\subseteq\F[z]^n$ be two submodules having the same rank and
satisfying $\hat{\cC}\subseteq\cC$.
Furthermore, let $\hat{\cC}$ be a \CC. Then $\hat{\cC}=\cC$.
\end{lemma}
\begin{proof}
Let $\cC=\im G$ and $\hat{\cC}=\im\hat{G}$ where $G,\,\hat{G}\in\F[z]^{k\times n}$
and~$\hat{G}$ is right invertible.
The assumption $\hat{\cC}\subseteq\cC$ implies the existence
of some matrix $U\in\F[z]^{k\times k}$ such that $\hat{G}=UG$. Using a right inverse
of $\hat{G}$ shows that $U\in Gl_k(\F[z])$ and the assertion follows.
\end{proof}

The complexity is also known as the {\em overall constraint length\/}~\cite[p.~55]{JoZi99},
\cite[p.~721]{Fo70}
or the {\em degree\/}~\cite[Def.~3.5]{McE98} of the code.
It is an important parameter describing the size of the code and of the encoding process.
In the coding literature a right invertible matrix is often also called
{\em basic\/}~\cite[p.~730]{Fo70} or
{\em delay-free and non-catastrophic}, see~\cite[p.1102]{McE98}.
Often in coding literature convolutional codes are defined as subspaces of
the vector space $\Flaurent^n$ of vector valued Laurent series over~$\F$, see for
instance~\cite{McE98} and~\cite{Fo70}.
However, as long as one restricts to right invertible generator matrices
it makes no difference with respect to code properties and code constructions
whether one works in the context of infinite message and codeword sequences
(Laurent series) or finite ones (polynomials).
Only for decoding it becomes important whether or not one may assume the
sent codeword to be finite.
The issue whether convolutional coding theory should be based
on finite or infinite message sequences, has first been raised and discussed in
detail in~\cite{RSY96,Ro01}.

It is well-known~\cite[Thm.~5]{Fo70} or~\cite[p.~495]{Fo75} that each submodule of
$\F[z]^n$ has a minimal generator matrix in the sense of the next definition.
In the same paper~\cite[Sec.~4]{Fo75} it has been shown how to derive such a matrix
from a given generator matrix in a constructive way.

\begin{defi}\label{D-minBasis}
\begin{arabiclist}
\item For $v=\sum_{j=0}^N v_jz^j\in\F[z]^n$ where $v_j\in\F^n$ and $v_N\not=0$
      let $\deg v:=N$ be the {\em degree\/} of~$v$.
      Moreover, put $\deg0=-\infty$.
\item Let $G\in\F[z]^{k\times n}$ be a matrix with rank~$k$ and complexity~$\delta$
      and let $\nu_1,\ldots,\nu_k$ be the degrees of the rows of~$G$.
      We say that $G$ is {\sl minimal\/} if $\delta=\sum_{i=1}^k\nu_i$.
      In this case, the row degrees of~$G$ are uniquely determined by the submodule
      $\cS:=\im G$. They are called the {\sl Forney indices\/}
      of~$\cS$.
\end{arabiclist}
\end{defi}
The notion ``minimal'' stems from the (simple) fact that for an arbitrary generator
matrix~$G$ one has $\delta\leq\sum_{i=1}^k \nu_i$.
Thus, in a minimal generator matrix the rows degrees have been reduced to their minimal
values.

From the above it follows that a convolutional code with parameters \nkd\ has a constant
generator matrix if and only if $\delta=0$.
In that case the code can be regarded as an $(n,k)$-block code.

The most important concept for a code is its distance.
It measures the error-correcting capability, hence the quality,  of the code.
The definition of the distance of a convolutional code is straightforward.
For a constant vector $w=(w_1,\ldots,w_n)\in\F^n$ we define, just like in block code theory,
its {\em (Hamming) weight\/} as $\wt(w)=\#\{i\mid w_i\not=0\}$.
For a polynomial vector $v=\sum_{j=0}^N v_j z^j\in\F[z]^n$, where $v_j\in\F^n$,
the {\em weight\/} is defined as $\wt(v)=\sum_{j=0}^N\wt(v_j)$.
Then the {\em (free) distance\/} of a code $\cC\subseteq\F[z]^n$ with generator matrix
$G\in\F[z]^{k\times n}$ is given as
\begin{equation}\label{e-dist}
   \dist(\cC):=\min\{\wt(v)\mid v\in\cC,\;v\not=0\}=
   \min\big\{\wt(uG)\,\big|\, u\in\F[z]^k,\;u\not=0\big\}.
\end{equation}
In coding theoretic terms, this notion is based on counting only the number of errors
during data transmission, but not their magnitude;
for more details about the distance of convolutional codes see for
instance~\cite[Sec.~3.1]{JoZi99}.
Although we will not present any theoretical results concerning the distance of a cyclic
convolutional code, we will show several examples of codes which do have optimal distance.
In all these cases the distances have been computed with a computer algebra program and then
compared to some suitable bound known from the literature.
One of these bound is the {\em generalized Singleton bound\/}~\cite{RoSm99} stating that the
distance~$d$ of a code with parameters \nkd\ over any field satisfies
\begin{equation}\label{e-MDS}
   d\leq S\nkd:=(n-k)\Big(\Big\lfloor\frac{\delta}{k}\Big\rfloor+1\Big)+\delta+1.
\end{equation}
A code~$\cC$ with $\dist(\cC)=S\nkd$ is called an MDS code.
The {\em Griesmer bound\/} also takes the field size into account.
It states that each code over a field with~$q$ elements and with parameters $(n,k,\delta)$ and
largest Forney index~$m$ has distance~$d$ bounded by
\begin{equation}\label{e-Griesmer}
  d\leq \max\Big\{d'\in\{1,\ldots,S\nkd\}\,\Big|\,\sum_{l=0}^{k(m+i)-\delta-1}
                    \Big\lceil\frac{d'}{q^l}\Big\rceil\leq n(m+i)
         \text{ for all }i\in\hat{\N}\Big\},
\end{equation}
see~\cite[3.22]{JoZi99} for $q=2$ and~\cite[Thm.~3.4]{GS03} for general field size.
Later we will present several codes where the distance attains this maximum value.

\section{Cyclic Convolutional Codes}
\setcounter{equation}{0}

In this section we will introduce the notion of cyclicity for \CC{}s.
After recalling from~\cite{Pi76} that the classical notion of invariance under cyclic shift
will always lead to complexity zero, we will introduce the skew-polynomial ring~$\Azs$,
isomorphic to $\F[z]^n$ as left $\F[z]$-module, and call the codes corresponding to
left ideals in~$\Azs$ cyclic.
We will briefly discuss some features of~$\Azs$ and summarize the main results about cyclic
codes, as obtained in~\cite{GS02}, in Theorem~\ref{T-ideals}.
From this we will derive that cyclic codes always have a cyclic direct complement, thereby
showing that the family of cyclic codes coincides with the family of those left ideals
in~$\Azs$ that are direct summands.

Just like for cyclic block codes we assume from now on that
\[
  \text{the length~$n$ and the field size~$|\F|$ are coprime.}
\]
Recall that a block code $\cC\subseteq\F^n$ is called cyclic if it is
invariant under the cyclic shift, i.~e.
\begin{equation}\label{e-cs}
  (v_0,\ldots,v_{n-1})\in\cC\Longrightarrow
  (v_{n-1},v_0,\ldots,v_{n-2})\in\cC
\end{equation}
for all $(v_0,\ldots,v_{n-1})\in\F^n$.
It is well-known that this is the case if and only if~$\cC$ is an ideal in
the quotient ring
\begin{equation}\label{e-A}
     A:=\F[x]/_{\DS \ideal{x^n-1}}
     =\Big\{\sum_{i=0}^{n-1}f_ix^i\;\mod(x^n-1)\,\Big|\,
          f_0,\ldots,f_{n-1}\in\F\Big\},
\end{equation}
canonically identified with $\F^n$ via
\[
  \p: \F^n\longrightarrow A,\quad
  (v_0,\ldots,v_{n-1})\longmapsto\sum_{i=0}^{n-1}v_ix^i\;\mod(x^n-1).
\]
Recall that the cyclic shift in~$\F^n$ translates
into multiplication by~$x$ in~$A$, i.~e.
\begin{equation}\label{e-cx}
  \p(v_{n-1},v_0,\ldots,v_{n-2})=x\p(v_0,\ldots,v_{n-1})
\end{equation}
for all $(v_0,\ldots,v_{n-1})\in\F^n$.
It is well-known that each ideal $I\subseteq A$ is principal,
hence there exists some $g\in A$ such that $I=\ideal{g}$.
One can even choose~$g$ as a monic divisor of $x^n-1$, in which case it is
usually called the {\em generator polynomial\/} of the code
$\p^{-1}(I)\subseteq\F^n$.

In order to extend the situation of cyclic block codes to the convolutional setting,
we have to replace the vector space~$\F^n$ by the free module $\F[z]^n$ and,
consequently, the ring~$A$ by the polynomial ring~$A[z]$ over~$A$.
Then we can extend the map~$\p$ above coefficient-wise to polynomials, thus
\begin{equation}\label{e-p}
  \p: \F[z]^n\longrightarrow A[z],\quad
  \sum_{j=0}^N z^jv_j\longmapsto \sum_{j=0}^N z^j\p(v_j),
\end{equation}
where, of course, $v_j\in\F^n$ and thus $\p(v_j)\in A$ for all~$j$.
This map is an isomorphism of $\F[z]$-modules.
Its inverse will be denoted by
\begin{equation}\label{e-v}
     \v:=\p^{-1}.
\end{equation}
Again, by construction the cyclic shift in $\F[z]^n$ corresponds to multiplication
by~$x$ in $A[z]$, that is, we have\eqnref{e-cx} for all
$(v_0,\ldots,v_{n-1})\in\F[z]^n$.
At this point it is quite natural to call a \CC\ $\cC\subseteq\F[z]^n$ cyclic if it is
invariant under the cyclic shift, i.~e. if\eqnref{e-cs} holds true for all
$(v_0,\ldots,v_{n-1})\in\F[z]^n$.
This, however, does not result in any codes other than block codes due to
the following result, see~\cite[Thm.~3.12]{Pi76} and \cite[Thm.~6]{Ro79}.
An elementary proof can be found at~\cite[Prop.~2.7]{GS02}.

\begin{theo}\label{T-CCCclassic}
Let $\cC\subseteq\F[z]^n$ be a convolutional code with parameters $(n,k,\delta)$
such that\eqnref{e-cs} holds true for all $(v_0,\ldots,v_{n-1})\in\F[z]^n$.
Then $\delta=0$, hence~$\cC$ is a block code.
\end{theo}

This result has led Piret~\cite{Pi76} to suggesting a different notion of
cyclicity for \CC{}s.
We will present this notion in the slightly more general version
introduced by Roos~\cite{Ro79}.

In order to do so notice that~$\F$ can be regarded as a subfield of the
ring~$A$ in a natural way. As a consequence,~$A$ is an $\F$-algebra.
In the sequel the automorphism group $\AutF(A)$ of the $\F$-algebra~$A$
will play an important role.
It is clear that each automorphism $\sigma\in\AutF(A)$ is uniquely
determined by the single value $\sigma(x)\in A$.
In particular, $\sigma(x)=x$ determines the identity map on~$A$.
But, of course, not every choice for $\sigma(x)$ determines an automorphism on~$A$.
Since~$x$ generates the $\F$-algebra~$A$, the same has to be true
for~$\sigma(x)$ and, more precisely, we obtain for $a\in A$
that $\sigma(x)=a$ determines an automorphism on~$A$ if and only if
$1,\,a,\ldots,\,a^{n-1}$ are linearly independent over~$\F$ and $a^n=1$.
A better way to determine $\AutF(A)$ will be described below in Remark~\ref{R-Aut}.

The main idea of Piret was to impose a new ring structure on $A[z]$ and to call a
code cyclic if it is a left ideal with respect to that ring structure.
The new structure is non-commutative and based on an (arbitrarily chosen)
automorphism on~$A$.
In detail, this looks as follows.

\newpage
\begin{defi}\label{D-CCC}
Let $\sigma\in\AutF(A)$.
\begin{arabiclist}
\item On the set $A[z]$ we define addition as usual while multiplication is defined
      via the rule
      \[
         \sum_{j=0}^Nz^j a_j\cdot\sum_{l=0}^M z^l b_l=
         \sum_{t=0}^{N+M}z^t\sum_{j+l=t}\sigma^l(a_j)b_l
         \text{ for all }N,\,M\in\N_0\text{ and }a_j,\,b_l\in A
      \]
      along with classical multiplication for the coefficients in the quotient ring~$A$.
      This turns $A[z]$ into a skew-polynomial ring, denoted by $\Azs$.
      We also call $\Azs$ a {\em Piret-algebra}.
\item Consider the map $\p:\F[z]^n\!\rightarrow\!\Azs$ as in\eqnref{e-p}, where now
      the images $\p(v)\!=\!\sum_{j=0}^Nz^j\p(v_j)$ are regarded as elements of
      $\Azs$.
      A submodule $\cS\subseteq\F[z]^n$ is said to be
      $\sigma$-{\em cyclic\/} if $\p(\cS)$ is a left ideal in $\Azs$.
      A \CC\ $\cC\subseteq\F[z]^n$ is said to be
      $\sigma$-{\em cyclic\/} if~$\cC$ is a direct summand of~$\F[z]^n$ and a \scy\
      submodule.
\end{arabiclist}
\end{defi}
A few comments are in order. First of all, notice that multiplication is
determined by the rule
\begin{equation}\label{e-az}
   az=z\sigma(a)\text{ for all }a\in A
\end{equation}
along with the rules of a non-commutative ring.
Hence, unless~$\sigma$ is the identity, the indeterminate~$z$ does not commute
with its coefficients.
Consequently, it becomes important to distinguish between left and right
coefficients of~$z$.
Of course, the coefficients can be moved to either side by applying the
rule\eqnref{e-az} since~$\sigma$ is invertible.
In the sequel we will always use the representation via right coefficients since
that is the one needed for the map~$\p$ in part~(2) above.
Since multiplication inside~$A$ remains the same as before~$A$ is a
commutative subring of $\Azs$. Moreover, since
$\sigma|_{\F}=\text{id}_{\F}$, the classical polynomial ring $\F[z]$ is a
commutative subring of $\Azs$, too.
As a consequence, $\Azs$ is a left and right $\F[z]$-module and the
map~$\p:\F[z]^n\rightarrow\Azs$ is an isomorphism of left $\F[z]$-modules
(but not of right $\F[z]$-modules).
In the special case where $\sigma=\text{id}_{A}$, the ring $\Azs$ is
the classical commutative polynomial ring and we know from
Theorem~\ref{T-CCCclassic} that no \scy\ \CC{}s with nonzero complexity
exist.
Finally, it should be noted that cyclic block codes (in the classical sense
of\eqnref{e-cs}) are \scy\ for all automorphisms~$\sigma$.

It is also worth being noted that, due to the definition above, $\sigma$-cyclic
convolutional codes are the left $\Azs$-submodules of $\Azs$ that are at the
same time direct summands of the left $\F[z]$-module $\Azs$.
As it will turn out this implies that they are direct summands
as $\Azs$-modules.
In other words, each $\sigma$-cyclic code has a direct complement that is
$\sigma$-cyclic, too (see~Corollary~\ref{C-dirSummands} below).

\begin{exa}\label{E-binarylength7}
Let us consider the case where $\F=\F_2$ and $n=7$. Thus
$A=\F_2[x]/_{\ideal{x^7-1}}$.
In this case $\AutF(A)$ contains~$18$ automorphisms (see also~\cite[p.~680, Table~II]{Ro79}),
one of which is defined via $\sigma(x)=x^5$.
We choose this automorphism for the following computations.
Consider the polynomial
\begin{equation}\label{e-polyg}
  g := 1+x^2+x^3+x^4+z(x+x^2+x^3+x^5)+z^2(1+x+x^4+x^6)\in\Azs
\end{equation}
and denote by $\lideal{g}:=\{fg\mid f\in\Azs\}$ the left ideal generated
by~$g$ in $\Azs$.
Moreover, put $\cC:=\v(\lideal{g})\subseteq\F[z]^7$.
We will show now that $\cC$ is a direct summand of $\F[z]^7$, hence
$\cC$ is a  \scy\ \CC.
In order to do so we first notice that
\[
   \lideal{g}=\spann_{\F[z]}\big\{g,\,xg,\ldots,x^6g\big\}
\]
and therefore, using the isomorphism~$\v$ from\eqnref{e-v},
\[
  \cC=\big\{u M\,\big|\, u\in\F[z]^7\big\}
  \text{ where }
       M=\begin{bmatrix}\v(g)\\\v(xg)\\\vdots\\\v(x^6g)
       \end{bmatrix}\in\F[z]^{7\times7}.
\]
Thus we have to compute $x^ig$ for $i=1,\ldots,6$.
Using the multiplication rule in\eqnref{e-az} we obtain
\begin{align*}
  xg&=x+x^3+x^4+x^5+z(1+x+x^3+x^6)+z^2(x+x^3+x^4+x^5),
  \\
  x^2g&=x^2+x^4+x^5+x^6+z(x+x^4+x^5+x^6)+z^2(1+x+x^2+x^5),
  \\
   x^3g&=1+x^3+x^5+x^6+z(x^2+x^3+x^4+x^6)+z^2(1+x^3+x^5+x^6)
        =g+x^2g.
\end{align*}
Since $x^3g$ is in the $\F$-span of the previous elements, we obtain
$\lideal{g}=\spann_{\F[z]}\big\{g,xg,x^2g\big\}$ and,
since~$\v$ is an isomorphism, $\cC=\big\{u G\,\big|\, u\in\F[z]^3\big\}$,
where
\[
  G=\begin{bmatrix}\v(g)\\\v(xg)\\\v(x^2g)
    \end{bmatrix}
   =\begin{bmatrix} 1+z^2&z+z^2&1+z&1+z&1+z^2&z&z^2\\
                    z&1+z+z^2&0&1+z+z^2&1+z^2&1+z^2&z\\
                    z^2&z+z^2&1+z^2&0&1+z&1+z+z^2&1+z
    \end{bmatrix}.
\]
One can easily check that the matrix~$G$ is right invertible and minimal
(see Definition~\ref{D-minBasis}).
Hence $\cC\subseteq\F[z]^7$ is indeed a \CCC.
It is worth mentioning that $\dist(\cC)=12$ (derived by a computer algebra program)
and this is the optimum value for any
convolutional code over~$\F_2$ with parameters $(7,3,6)$ by virtue of the
Griesmer bound\eqnref{e-Griesmer}.
\end{exa}

In order to proceed with the theory of \CCC s one needs some knowledge about the left
ideals in the skew-polynomial ring $\Azs$.
In particular, we need to understand whether a given left ideal corresponds to a
convolutional code rather than just to a submodule and, if so, if the parameters
(dimension and complexity) can be recovered from the ideal.
All this has been answered in the affirmative in~\cite{GS02}.
In the sequel we will present the according results.

The main tool for describing the left ideals in $\Azs$ is the fact that~$A$
is a semi-simple ring.
Since we need the details of this fact we will first elaborate on this.
By comprimeness of the length~$n$ and the field size~$|\F|$, the polynomial $x^n-1$ is square free,
say
\begin{equation}\label{e-xn-1}
   x^n-1=\pi_1\cdot\ldots\cdot \pi_r,
\end{equation}
where $\pi_1,\ldots,\pi_r\in\F[x]$ are irreducible, monic, and pairwise different.
We will also assume that the polynomials are ordered according to
\begin{equation}\label{e-ordering}
  \deg_x\pi_1=\ldots=\deg_x\pi_{r_1}<\ldots<\deg_x\pi_{r_1+\ldots+r_{s-1}+1}
   =\ldots=\deg\pi_{r_1+\ldots+r_s},
\end{equation}
where $r_1+\ldots+r_s=r$.
Using $r_0:=0$ and $l_t:=\sum_{\lambda=0}^{t-1}r_{\lambda}+1$ for $t=1,\ldots,s$,
we have the partition
$\{1,\ldots,r\}=R^{(1)}\cup\ldots\cup R^{(s)}$ where
$R^{(t)}=\{l_t,\,l_t+1,\ldots,l_t+r_t-1\}$.
It will also be convenient to use equivalence relation
\begin{equation}\label{e-equivind}
    k\equiv l :\Longleftrightarrow \deg_x\pi_k=\deg_x\pi_l.
\end{equation}
Hence $k\equiv l$ if and only if~$k$ and~$l$ belong to the same index set
$R^{(t)}$ for some~$t$.

The Chinese Remainder Theorem provides us with an isomorphism of rings
\begin{equation}\label{e-CRT}
  \psi: A\longrightarrow K_1\times\ldots\times K_r,\quad
         a\longmapsto \db{\rho_1(a),\ldots,\rho_r(a)},
\end{equation}
where $K_k=\F[x]/_{\DS\ideal{\pi_k}}$ and $\rho_k$ denotes the canonical projection.
Notice that $K_k\cong K_l$ if and only if $k\equiv l$.
As indicated in\eqnref{e-CRT}, the elements in the direct product will be denoted by
$\db{a_1,\ldots,a_r}$.
It is easy to see that the elements
\[
     \ve{k}:=\psi^{-1}\big(\db{(\delta_{kj})_{1\leq j\leq r}}\big) \text{ for }
     k=1,\ldots,r
\]
form the uniquely determined set of primitive idempotents in~$A$.
We call the subfield
$K^{(k)}:=\ve{k}A=\psi^{-1}(0\times\ldots\times0\times K_k\times0\times\ldots\times0)$
the $k$-th component of~$A$.
Obviously, $A=K^{(1)}\oplus\ldots\oplus K^{(r)}$, showing that~$A$ is a semisimple
left-Artinian ring, see e.~g.~\cite[Ch.~IX, Sec.~3.1]{Hu74}.
In particular, $A$ has only finitely many ideals, each of which being isomorphic to a
direct product of fields.
Moreover,
\begin{equation}\label{e-unitsA}
   a\in A \text{ is a unit in }A \Longleftrightarrow
   \ve{l}a\not=0\text{ for all }l=1,\ldots,r.
\end{equation}

Let us now study the effect of a given automorphism $\sigma\in\AutF(A)$ on the
components.
It is straightforward to see that for each~$k$ we have $\sigma(K^{(k)})=K^{(l)}$
for some~$l$ such that $l\equiv k$.
In other words,
\begin{equation}\label{e-sigmaeps}
   \sigma(\ve{k})=\ve{l} \text{ for some $l$ such that $k\equiv l$.}
\end{equation}
This gives rise to the following definition.

\begin{defi}\label{D-sigmapermut}
Let $\sigma\in\AutF(A)$.
Define the permutation $\Pi_{\sigma}\in S_r$ via
$\Pi_{\sigma}(k)=l$ where~$l$ is such that $\sigma(\ve{k})=\ve{l}$ for all
$k=1\,\ldots,r$.
We call $\Pi_{\sigma}$ the permutation corresponding to~$\sigma$.
Furthermore, define the equivalence relation $\equivs$ on the index set
$\{1,\ldots,r\}$ via $k\equivs l$ if there exists some $i\in\N_0$ such that
$\sigma^i(\ve{k})=\ve{l}$.
\end{defi}
Of course, the permutation $\Pi_{\sigma}$ simply reflects the permutation induced
by~$\sigma$ on the set $\{\ve{1},\ldots,\ve{r}\}$, that is,
$\sigma(\ve{k})=\ve{\Pi_{\sigma}(k)}$.
The equivalence relation $\equivs$ can also be expressed as $k\equivs l$ if and only if $k$
and~$l$ belong to the same cycle of the permutation~$\Pi_{\sigma}$.
Since the permutation~$\Pi_{\sigma}$ satisfies
$\Pi_{\sigma}(R^{(t)})=R^{(t)}$ for all $t=1,\ldots,r$, see\eqnref{e-sigmaeps}, we obtain that
each of its cycles is contained in one of the sets $R^{(t)}$.
In other words
\[
     k\equivs l\Longrightarrow k\equiv l\quad \text{ for all }k,\,l\in\{1,\ldots,r\}.
\]

The consideration above provides us with an alternative way to compute the
automorphisms on~$A$.

\begin{rem}\label{R-Aut}
It is straightforward to see that each permutation $\Pi\in S_r$ satisfying
$\Pi(R^{(k)})=R^{(k)}$ for all~$k\in\{1,\ldots,r\}$ is the
permutation~$\Pi_{\sigma}$ of an $\F$-automorphism~$\sigma$ on~$A$.
Hence~$\sigma$ is such that $\sigma(K^{(k)})=K^{(\Pi(k))}$ for all $k=1,\ldots,r$.
Since there are, in general, many isomorphisms between~$K^{(k)}$ and~$K^{(\Pi(k))}$,
the permutation~$\Pi$ does not completely determine the automorphism.
Rather, we obtain all automorphisms~$\sigma$ on~$A$ satisfying $\Pi_{\sigma}=\Pi$
by fixing one isomorphism between $K^{(k)}$ and $K^{(\Pi(k))}$ and using the automorphism
group $\AutF(K^{(k)})$ for presenting the remaining ones.
One can show that in this way one obtains all automorphisms on~$A$, see~\cite{Ve85}.
With this consideration one can easily compute the cardinality of the automorphism
group.
Indeed, notice that $r_1!\cdots r_s!$ counts the number of all permutations~$\Pi$ satisfying
$\Pi(R^{(t)})=R^{(t)}$ for all~$t$.
Since each~$k$ is in one of the sets $R^{(t)}=\{l_t,l_t+1,\ldots,l_t+r_t-1\}$ and
$|\AutF(K^{(l_t)})|=\deg_x\pi_{l_t}$ the above leads to
$|\AutF(A)|=
  (\deg_x\pi_{l_1})^{r_1}\cdots(\deg_x\pi_{l_s})^{r_s}r_1!\cdots r_s!$.
For more details see~\cite[Sec.~3]{GS02}.
\end{rem}

Having this description of the semi-simple ring~$A$ and its automorphisms available
we will now fix some $\sigma\in\AutF(A)$ and turn to the Piret-algebra $\Azs$ over~$A$.
This ring is, of course, an~$A$-module and as such semisimple
(i.~e. every $A$-submodule of $\Azs$ is a direct summand),
see~\cite[Ch.~IX, Thm.~3.7]{Hu74}.
However, for our investigation of \scy\ codes we need to understand the ring structure
along with the left $\F[z]$-module structure.
This has been worked out in detail in~\cite{GS02} and leads to the following.

Using $1=\ve{1}+\ldots+\ve{r}$ we can write each polynomial $f\in\Azs$ in the form
\[
   f=f^{(1)}+\ldots+f^{(r)},
   \text{ where }f^{(k)}:=\ve{k}f.
\]
We call $f^{(k)}$ the {\em $k$-th component\/} of~$f$. Furthermore, the set
$T_f:=\big\{k\in\{1,\ldots,r\}\,\big|\, f^{(k)}\not=0\big\}$ is called the
{\em support\/} of $f$.
From\eqnref{e-az} it follows that $\ve{k}z^{\mu}=z^{\mu}\ve{k'}$
for some $k'$ such that $k\equivs k'$.
Therefore, each $f\in\Azs$ can be written as an $A$-linear combination of the
elements
\begin{equation}\label{e-monom}
    z^{\mu}\ve{k},\ \mu\geq0,\ k=1,\ldots,r.
\end{equation}
We call these elements the {\em monomials\/} of $\Azs$.
In particular, the $k$-th component $f^{(k)}=\ve{k}f$ of~$f$ satisfies
\begin{equation}\label{e-fkcoeff}
  f^{(k)}\in\spann_{A}\{ z^{\mu}\ve{k'}\mid \mu\geq0,\ k'\equivs k\}
\end{equation}
(where the span has to be understood with respect to right coefficients).
Thus, the (right) coefficients of $f^{(k)}$ are not in $\ve{k}A$ but rather move
around in the fields $K^{(k')}=\ve{k'}A$, where $k'\equivs k$.
From this and the orthogonality of the idempotents it follows immediately
the orthogonality of components corresponding to disjoint cycles, precisely
\begin{equation}\label{e-orthog}
  f,\,g\in\Azs,\ k\nequivs l\Longrightarrow f^{(k)}g^{(l)}=g^{(l)}f^{(k)}=0.
\end{equation}

\begin{exa}\label{E-binarylength7a}
Consider again Example~\ref{E-binarylength7} where $\F=\F_2,\,n=7$ and $\sigma(x)=x^5$.
The polynomial $x^7-1$ decomposes into $x^7-1=\pi_1\pi_2\pi_3$ where
\[
   \pi_1=x+1,\ \pi_2=x^3+x+1,\ \pi_3=x^3+x^2+1.
\]
Thus, in the notation of\eqnref{e-xn-1} and\eqnref{e-ordering}, $r=3,\,s=2$ and
$R^{(1)}=\{1\}$ and $R^{(2)}=\{2,\,3\}$.
Furthermore, one has the primitive idempotents
\[
  \ve{1}=1+x+x^2+x^3+x^4+x^5+x^6,\ \ve{2}=1+x+x^2+x^4,\ \ve{3}=1+x^3+x^5+x^6,
\]
which can easily be checked by verifying that
$(\ve{k}\,\mod\pi_i)=\delta_{ik}$ for $i,\,k=1,\,2,\,3$.
Moreover,
$\sigma(\ve{1})=\ve{1},\,\sigma(\ve{2})=\ve{3},\,\sigma(\ve{3})=\ve{2}$. In other
words,~$\sigma$ induces the permutation $\Pi_{\sigma}=(1)(2,3)$.
It can straightforwardly be shown that the polynomial~$g$ given in\eqnref{e-polyg}
satisfies $g^{(1)}=0=g^{(2)}$ as well as
\[
  g=g^{(3)}=
  \ve{3}(1+x+x^2)+z\ve{2}x+z^2\ve{3}x.
\]
Hence $\psi(g)=\db{0,0,1+x+x^2}+z\db{0,x,0}+z^2\db{0,0,x}$.
This can be verified directly and expresses the fact that the
coefficient~$g_0$ of~$z^0$ in~$g$ satisfies $(g_0\,\mod\pi_1)=0=(g_0\,\mod\pi_2)$
and~$(g_0\,\mod\pi_3)=1+x+x^2$.
According relations hold for the coefficients of $z$ and~$z^2$.
\end{exa}

Having this description of the polynomials in the Piret-algebra $\Azs$ at hand we are now
in a position to investigate the left ideals.
In~\cite{GS02} a Groebner-type theory has been established for $\Azs$.
It is based on the monomials given in\eqnref{e-monom} and leads to a reduction
algorithm just like for commutative polynomials in several variables.
This looks as follows.

\begin{defi}\label{D-monomorder}
\begin{alphalist}
\item Given two monomials $z^{\mu}\ve{k}$ and $z^{\nu}\ve{l}$ we define
      \[
         z^{\mu}\ve{k}<z^{\nu}\ve{l}\Longleftrightarrow
         \mu<\nu \text{ or }\mu=\nu\text{ and } k<l.
      \]
\item For a polynomial
      $f=\sum_{\nu\geq0}z^{\nu}f_{\nu}=\sum_{\nu\geq0}\sum_{l=1}^rz^{\nu}\ve{l}f_{\nu}\in\Azs$
      define $LM(f)$ to be the
      largest monomial $z^{\mu}\ve{k}$ (with respect to $<$) which has a nonzero coefficient
      in~$f$, that is, for which $\ve{k}f_{\mu}\not=0$.
      We call $LM(f)$  the {\em leading monomial of~$f$}.
      The summands $z^{\nu}\ve{l}f_{\nu}$ are called the {\em terms of~$f$}.
\item A polynomial $f\in\Azs$ is called {\em (left) reduced\/} if for all $k,\,l=1,\ldots,r$,
      where $k\not=l$, no nonzero term of $f^{(k)}$ is right divisible by $LM(f^{(l)})$.
\item A polynomial $f\in\Azs$ is called a {\em component\/} if $f=f^{(k)}$ for some
      $k=1,\ldots,r$.
\end{alphalist}
\end{defi}
One easily verifies that $<$ is a well-ordering on the set of monomials with respect
to multiplication as far as the result is nonzero.
Notice that a component $f^{(k)}$ is always reduced.

In~\cite{GS02} a reduction procedure for polynomials has been
established which, just like in the commutative case of several variables, leads in a
constructive way to a type of Groebner bases for left ideals in $\Azs$.
We will need the following results on principal left ideals.
They have been proven in~\cite[Thm.~4.5, Cor.~4.13(b), Prop.~7.10, Thm.~7.13]{GS02}.

\begin{theo} \label{T-ideals}
Fix $\sigma\in\AutF(A)$. Then
\begin{arabiclist}
\item Each principal left ideal $\cI\in\Azs$ has a reduced generator polynomial.
      Precisely, there exists a reduced polynomial $g\in\Azs$ such that
      \[
          \cI=\lideal{g}:=\{fg\mid f\in\Azs\}.
      \]
      Moreover, the reduced generator is unique up to left multiplication by units
      in~$A$.
\item Let $\cC\subseteq\F[z]^n$ be a \scy{} \CC.
      Then the associated left ideal~$\p(\cC)$ is principal and thus has a reduced
      generator $g\in\Azs$. Moreover, the support of~$g$ satisfies $T_g=T_{g_0}$
      where~$g_0$ denotes the constant term of~$g$.
\item Let $g\in\Azs$ be a reduced polynomial.
      Then $\v(\lideal{g})\subseteq\F[z]^n$ is a direct summand of $\F[z]^n$
      (thus a \scy{} \CC) if and only if there exist $a\in A$ and a unit $v\in\Azs$
      such that $g=av$.
\item Let $g\in\Azs$ be a reduced polynomial with support~$T_g$.
      For $l\in T_g$ let $\deg_x\pi_l=\kappa_l$, where $\pi_l$ is as in\eqnref{e-xn-1},
      and put $\kappa:=\sum_{l\in T_g}\kappa_l$.
      Then the matrix
      \begin{equation}\label{e-directsum}
          G:=\begin{bmatrix} \v\big(x^ig^{(l)}\big)\end{bmatrix}_{l\in
          T_g,\,i=0,\ldots,\kappa_l-1}\in\F[z]^{\kappa\times n}
      \end{equation}
      is a minimal generator matrix of the submodule~$\cS:=\im G\subseteq\F[z]^n$.
      As a consequence,~$\cS$ is a submodule of rank~$\kappa$ and complexity
      $\delta=\sum_{l\in T_g}\kappa_l \deg_z g^{(l)}$.
      The Forney-indices are given by the numbers $\deg_z g^{(l)},\;l\in T_g$,
      each one counted $\kappa_l$ times.
\end{arabiclist}
\end{theo}
We wish to comment on these results.
First of all, it is worth mentioning that $\Azs$ is not a left principal ideal ring.
Part~(2) above only states that left ideals associated to direct summands in
$\F[z]^n$ are principal.
Indeed, there exist left ideals that are not principal~\cite[Exa.~4.6(a)]{GS02}.
Secondly, as for part~(3) above we wish to mention that each left inverse of some
$v\in\Azs$ is also a right inverse \cite[p.~32]{GS02}; this will slightly simplify the
investigation of units.
Due to zero divisors in the coefficient ring~$A$, the skew-polynomial ring has
plenty of units of higher $z$-degree, i.~e., units, that are not in~$A$.
We will investigate this issue in more detail in the next section.
Notice that a unit itself is never reduced unless it is in~$A$, i.~e., a constant.
This follows for instance from~(1) since a unit generates (as a left ideal) the full
Piret-algebra $\Azs$, which in turn has the reduced polynomial~$1\in A$
as a generator.
Finally, we want to emphasize that according to~(4) the parameters of \scy\ \CC{}s
can occur only in certain combinations.
In particular, the Forney indices appear, in general, with higher multiplicities
depending on the degrees of the prime factors~$\pi_l$.
In the next section we will investigate this situation in more detail.

It is worth being stressed that part~(2) and~(3) above deal with direct summands of
the left module $\Azs$ over the ring~$\F[z]$ and not over $\Azs$.
However, it can easily be deduced from the above that direct summands with respect to
these different structures coincide. Indeed,

\begin{cor}\label{C-dirSummands}
Let $\cI$ be a left ideal in $\Azs$. Then the following are equivalent
\begin{romanlist}
\item $\cI$ is a direct summand of the left $\F[z]$-module $\Azs$,
\item $\cI$ is a direct summand of the left $\Azs$-module $\Azs$.
\end{romanlist}
In particular, a $\sigma$-cyclic code has a direct summand that is $\sigma$-cyclic
again.
Furthermore, if~$\cI$ is a direct summand, then $\cI=\lideal{g}$ where
$g=\sum_{l\in T_g}u^{(l)}$ for some unit $u\in\Azs$.
In this case, a direct complement is given by $\lideal{g'}$ where
$g':=\sum_{l\not\in T_g}u^{(l)}$.
\end{cor}
Before we give the proof we wish to add that, using the reduction procedure established
in~\cite[Sec.~4]{GS02}, it is possible to test
constructively whether or not a given reduced polynomial generates a direct summand.
This also produces a direct summand in a constructive way.
\\
\begin{proof}
The direction (ii)~$\Rightarrow$~(i) is clear since each $\Azs$-module is also an
$\F[z]$-module.
\\
(i)~$\Rightarrow$~(ii):
By Theorem~\ref{T-ideals}(2) the ideal~$\cI$ is principal, say
$\cI=\lideal{\hat{g}}$, where~$\hat{g}$ is a reduced polynomial.
By part~(3) of that theorem we have $\hat{g}=au$ for some $a\in A$ and
a unit $u\in\Azs$.
Then $T_a=T_{\hat{g}}$.
We can normalize the factor~$a$ in the following way.
Since the ring~$A$ is a direct product of fields, there exists a unit
$\hat{a}\in A$ such that $\hat{a}a=\sum_{l\in T_a}\ve{l}$.
Hence $\cI=\lideal{g}$ where $g:=\hat{a}au=\sum_{l\in T_a}u^{(l)}$ and
$T_g=T_a$.
A direct complement is given by the left ideal $\lideal{g'}$ where
$g':=\sum_{l\not\in T_g}u^{(l)}$.
In order to see this, notice first that $g+g'=u$ is a unit in $\Azs$ and hence
$\lideal{g}+\lideal{g'}=\Azs$.
Suppose now that $fg=f'g'\in\lideal{g}\cap\lideal{g'}$ for some $f,\,f'\in\Azs$.
Then $(f\sum_{l\in T_g}\ve{l}-f'\sum_{l\not\in T_g}\ve{l})u=0$ and,
since~$u$ is a unit, $f\sum_{l\in T_g}\ve{l}=f'\sum_{l\not\in T_g}\ve{l}$.
But this implies
$f\ve{k}=f\sum_{l\in T_g}\ve{l}\ve{k}=f'\sum_{l\not\in T_g}\ve{l}\ve{k}=0$
for all $k\in T_g$.
Hence $fg=\sum_{k\in T_g}f\ve{k}g=0$, showing that
$\lideal{g}\cap\lideal{g'}=\{0\}$.
All this also proves the additional assertion.
\end{proof}

The above shows that the set of $\sigma$-cyclic codes is the same as
the set of direct summands of the ring $\Azs$.
In this context it is worth being recalled that in every ring~$R$ with~$1$, a left
ideal~$\cI$, that is a direct summand (as left $R$-module), is left principal and even
has an idempotent generator.
We wish to emphasize that reduced generators, as guaranteed by
Theorem~\ref{T-ideals}(2), are in general not idempotent.
But the corollary above shows how idempotent generators can easily be obtained from the
reduced generator.
Indeed, with the data as in the corollary we have that $g+g'=u$ is a unit in $\Azs$.
Thus $1=u^{-1}g+u^{-1}g'$ and
$u^{-1}gu^{-1}g'=u^{-1}g-u^{-1}gu^{-1}g\in\lideal{g}\cap\lideal{g'}=\{0\}$.
From this it follows that both terms $u^{-1}g$ and $u^{-1}g'$ are idempotent generators of the
respective left ideal.
In general these idempotent generators have much higher degree than the reduced ones.
At any rate, as Theorem~\ref{T-ideals}(4) shows, the reduced generators are the more useful ones
when it comes to the associated module in $\F[z]^n$.

Since the reduced generator of a principal left ideal is essentially unique, the
following definition is well-posed.

\begin{defi}\label{D-supp}
Let $g\in\Azs$ be a reduced polynomial.
Then its support~$T_g$ is called the {\em support of the left ideal\/}~$\lideal{g}$
and also the {\em support of the submodule\/}~$\v(\lideal{g})$.
\end{defi}

The previous examples illustrate the results given so far.

\begin{exa}\label{E-binarylength7b}
Let us return once more to Example~\ref{E-binarylength7} and its
continuation in Example~\ref{E-binarylength7a}.
In that case the polynomial $g=g^{(3)}$
is reduced since it is a component. It generates a left ideal
corresponding to a code of rank~$3=\deg_x\pi_3$ and complexity
$6=\deg_x\pi_3\deg_z g^{(3)}$ which has been given explicitly in
Example~\ref{E-binarylength7}. This is compliant with what has
been stated in Theorem~\ref{T-ideals}(4). A $\sigma$-cyclic direct
complement of $\lideal{g}$ in $\Azs$ is given by the left ideal
generated by the polynomial
\[
  g'=x+x^3+x^4+z(1+x^3+x^5+x^6).
\]
One way to check this is by showing that $v=g+g'=1+x+x^2+z(1+x+x^2+x^6)+z^2(1+x+x^4+x^6)$
is a unit in $\Azs$.
This is indeed the case, its inverse is given by
$v^{-1}=1+x^2+x^3+x^6+z(x+x^2)+z^2(1+x^2+x^5+x^6)$.
The components of~$v$ are given by
\[
  v^{(1)}=\ve{1},\ v^{(2)}=\ve{2}(1+x+x^2)+z\ve{3},\
  v^{(3)}=\ve{3}(x^2+x+1)+z\ve{2}x+z^2\ve{3}x
\]
showing that $g=v^{(3)}$ while one easily verifies that $g'=v^{(1)}+v^{(2)}$.
We do not discuss how one obtains such a direct complement, since that needs
more detailed results from~\cite{GS02}.
\end{exa}

\section{Minimal Cyclic Codes}
\setcounter{equation}{0}
As before, let~$\F$ be a finite field such that~$n$ and~$|\F|$ are coprime and let
$\sigma\in\AutF(A)$ be a fixed automorphism, where~$A$ is as in\eqnref{e-A}.
In this section we will investigate the building blocks of \scy\ \CC{}s, the minimal
cyclic codes.
We will derive necessary and sufficient conditions for the automorphism~$\sigma$ to
allow for \scy\ codes with arbitrarily prescribed Forney indices.

As we saw in Theorem~\ref{T-ideals}(2) each \scy\ \CC\ $\cC\subseteq\F[z]^n$
corresponds to a principal left ideal in $\Azs$ which is generated by a reduced polynomial.
We will call each such reduced generator a {\em generator polynomial\/} of the
code~$\cC$.
Furthermore, part~(4) of that theorem shows that each \scy{} \CC\ can be presented
as the direct sum of \scy\ codes with components as generator polynomials.
Indeed, using the isomorphism~$\p$, Equation\eqnref{e-directsum} translates into the
direct sum
\[
    \lideal{g}=\bigoplus_{l\in T_g}\lideal{g^{(l)}}
\]
of left ideals in $\Azs$.
This leads to the following definition.

\begin{defi}\label{D-minC}
Let $\{0\}\not=\cC\subseteq\F[z]^n$ be a \scy{} \CC\ with generator polynomial $g\in\Azs$.
Then $\cC$ is called {\em minimal\/} if $g$ is a component, i.~e. if
$g=\ve{l}g=g^{(l)}$ for some $l\in\{1,\ldots,r\}$.
\end{defi}
The notion ``minimal'' (which is not related to minimal generator matrices) is justified
by the following result.

\begin{prop}\label{P-minC}
Let $\cC\subseteq\F[z]^n$ be a \scy{} \CC\ with generator polynomial $g\in\Azs$.
Then the following are equivalent.
\begin{romanlist}
\item $\cC$ is minimal,
\item $\cC\not=\{0\}$ and $\cC$ contains no proper \scy{} subcodes. Precisely,
      if $\hat{\cC}$ is a \scy{} \CC\ and $\{0\}\not=\hat{\cC}\subseteq\cC$, then
      $\hat{\cC}=\cC$.
\item There exists a unit $u\in\Azs$ such that $g=u^{(l)}$ for some index~$l$.
\end{romanlist}
\end{prop}
\begin{proof}
(i)~$\Rightarrow$~(ii):
By assumption $0\not=g=g^{(l)}$ for some index~$l$.
Let $\{0\}\not=\hat{\cC}$ be a \scy\ \CC\ with generator polynomial~$h\not=0$ and let
$\hat{\cC}\subseteq\cC$.
Then $\lideal{h}\subseteq\lideal{g}$, thus $h=fg$ for some $f\in\Azs$.
This implies $h_0=f_0g_0$ for the constant terms of the polynomials.
From Theorem~\ref{T-ideals}(2) we know that $g_0=g_0^{(l)}\not=0$, hence
$h_0=f_0\ve{l}g_0=h_0^{(l)}$.
Using again Theorem~\ref{T-ideals}(2) we deduce that $T_h=T_{h_0}=\{l\}$.
Thus $h=h^{(l)}$ and by Theorem~\ref{T-ideals}(4)
the codes $\hat{\cC}$ and $\cC$ have the same rank.
From Lemma~\ref{L-directsum} we conclude that $\hat{\cC}=\cC$.
\\
(ii)~$\Rightarrow$~(i):
follows directly from Theorem~\ref{T-ideals}(4) since each component of the
generator polynomial gives a \scy\ subcode of~$\cC$.
\\
The equivalence (i)~$\Leftrightarrow$~(iii) is clear with Corollary~\ref{C-dirSummands}.
\end{proof}

In the sequel we will show which parameters $(n,k,\delta)$ a minimal \scy\ \CC\
can attain.
From Theorem~\ref{T-ideals}(4) and Proposition~\ref{P-minC} we have the following
situation.

\begin{rem}\label{R-minCparam}
\begin{alphalist}
\item Any component~$u^{(l)}$ of a unit $u\in\Azs$ defines a minimal \scy\
      code $\v(\lideal{u^{(l)}})$ with parameters $(n,k,dk)$ where
      $k=\deg_x\pi_l$ and $d=\deg_z u^{(l)}$.
\item Any minimal \scy\ code in $\F[z]^n$ with support~$\{l\}$ has parameters
      $(n,k,dk)$ and Forney index~$d$ counted~$k$ times, where $k=\deg_x\pi_l$
      and~$d$ is the degree of the $l$-th component of a unit in $\Azs$.
\end{alphalist}
\end{rem}
Hence the question raised above amounts to investigating as to which degrees can occur
for a given component of a unit in $\Azs$.
The case where the complexity is zero is, of course, known from block code theory.
Indeed, for each $k\in\{\deg_x\pi_1,\ldots,\deg_x\pi_r\}$ there exists a cyclic
block code with parameters $(n,k)$, hence a \scy\ \CC\ with parameters $(n,k,0)$
for any automorphism~$\sigma$.
This follows also immediately from Remark~\ref{R-minCparam}(a).
The existence of \scy\ \CC{}s with nonzero complexity however, implies certain
relations between the parameters and the automorphism. Indeed, we have

\begin{lemma}\label{L-compl0}
Let $\cC\subseteq\F[z]^n$ be a minimal \scy\ code with generator polynomial
$g=g^{(l)}$.
Then~$\cC$ has complexity zero if and only if $g=g\ve{l}$.
Furthermore, if~$\cC$ has nonzero complexity then $\sigma\big(\ve{l}\big)\not=\ve{l}$.
\end{lemma}
\begin{proof}
If $\cC$ has complexity zero, then, by Theorem~\ref{T-ideals}(4), the polynomial~$g$ has
degree zero, thus $g\in A$.
But then $g=\ve{l}g=g\ve{l}$ follows from commutativity of~$A$.
Conversely, $g=\ve{l}g=g\ve{l}$ implies $\lideal{g}\subseteq\lideal{\ve{l}}$ and
thus $\cC\subseteq\v(\lideal{\ve{l}})$.
Both submodules are direct summands and by virtue of Theorem~\ref{T-ideals}(4)
they have the same rank.
Thus, Lemma~\ref{L-directsum} implies $\cC=\v(\lideal{\ve{l}})$ and therefore
has complexity zero.
As for the last assertion, notice that the identity $\sigma(\ve{l})=\ve{l}$ and the
very definition of multiplication in the Piret-algebra implies that~$\ve{l}$
is in the center of $\Azs$.
Hence $g=\ve{l}g=g\ve{l}$ and the code has complexity zero.
\end{proof}

As a consequence we have that for given parameters~$n$ and~$|\F|$ a given automorphism
$\sigma\in\AutF(A)$ admits (minimal) \scy\ \CC{}s of positive complexity only if the
permutation $\Pi_{\sigma}\in S_r$ is nontrivial.
This in turn is possible only if at least one of the sets $R^{(t)}$ contains more
than one element (see Definition~\ref{D-sigmapermut}) or in other words, if
$x^n-1$ has (at least) two prime factors of the same degree.
Recall that one easily obtains the degrees of the prime factors of $x^n-1$ by
computing the cyclotomic cosets modulo~$n$ over~$\F$, see~\cite[Ch.~7, \S~5]{MS77}.
With different methods it has been shown in~\cite[Sec.~VI]{Ro79} and
in~\cite[Prop.~3.4]{GS02} that the condition $\Pi_{\sigma}\not=\id$ is not only
necessary but also sufficient for the existence of \scy\ codes with positive complexity.
Our goal is to prove even more.
We will show that for any $\sigma\in\AutF(A)$ and any $l\in\{1,\ldots,r\}$ such that
$\sigma(\ve{l})\not=\ve{l}$ and for any~$d\in\N$ there exists a minimal \scy\ code
with parameters $(n,k,kd)$ where $k=\deg_x\pi_l$.
To this aim we need

\begin{defi}\label{D-lorder}
Let $\sigma\in\AutF(A)$ and $l\in\{1,\ldots,r\}$.
We define the $l$-order of~$\sigma$ as
$o_l(\sigma):=\min\{m\in\N\mid \sigma^m(\ve{l})=\ve{l}\}$.
\end{defi}
Using the permutation $\Pi_{\sigma}\in S_r$ associated with~$\sigma$, the $l$-order
can also be expressed as $o_l(\sigma)=\min\{m\in\N\mid \Pi_{\sigma}^m(l)=l\}$.
In other words, the $l$-order of $\sigma$ is the length of the cycle of
$\Pi_{\sigma}$ containing~$l$; therefore
\begin{equation}\label{e-lorder}
  l\equivs l'\Longrightarrow o_l(\sigma)=o_{l'}(\sigma).
\end{equation}

With the following lemma we will establish the existence of certain simple
units in $\Azs$.
They will suffice to show the existence of the desired minimal \scy\ codes.
We will also obtain that each unit in $\Azs$ can be expressed as a
finite product of these simple units.
In this sense we can construct, at least theoretically, all units of~$\Azs$ and
thus, by Corollary~\ref{C-dirSummands}, all \scy\ \CC{}s.

\newpage
\begin{lemma}\label{L-eltunits}
Let $\sigma\in\AutF(A)$ with $l$-order $o_l:=o_l(\sigma)$ where
$l\in\{1,\ldots,r\}$.
\begin{alphalist}
\item Let $a\in A$ and $d\in\N$.
      Put $u_{d,a,l}:=1+z^da\ve{l}\in\Azs$. Then
      \[
         u_{d,a,l}\text{ is a unit in }\Azs
         \Longleftrightarrow
         \left\{\begin{array}{ll}
         a^{(l)}\not=-\ve{l},&\text{ if }d=0,\\
         a^{(l)}=0 \text{ or }o_l\nmid d,&\text{ if }d>0.
         \end{array}\right.
      \]
      If $u_{d,a,l}$ is a unit in $\Azs$, then its inverse is given by $u_{d,-a,l}$.
      In this case we call $u_{d,a,l}$ an elementary unit.
\item Any unit in $\Azs$ can be written as a finite product of elementary units.
\end{alphalist}
\end{lemma}
\begin{proof}
(a) If $d=0$ then $u_{d,a,l}=1+a^{(l)}$ and the assertion follows from\eqnref{e-unitsA}.
Thus let $d>0$. We may assume $a^{(l)}\not=0$ for otherwise the assertion is trivial.
\\
``$\Rightarrow$'' Write $u:=u_{d,a,l}$, for short.
Since~$u$ is a unit, we know from Remark~\ref{R-minCparam}(a)
that $\v(\lideal{u^{(l)}})$ is a minimal \scy\ convolutional code
and its complexity is given by $\deg_x\pi_l\deg_z u^{(l)}$.
If $o_l\mid d$, then $\ve{l}z^d=z^d\ve{l}$ and thus
$u^{(l)}=\ve{l}u=u\ve{l}=\ve{l}+z^da^{(l)}$, hence $\deg_zu^{(l)}=d>0$.
But on the other side Lemma~\ref{L-compl0} implies that the complexity
of $\v(\lideal{u^{(l)}})$ is zero, a contradiction.
\\
``$\Leftarrow$'' Let $o_l\nmid d$.
Then $\sigma^d(\ve{l})\not=\ve{l}$ and thus $\sigma^d(\ve{l})\ve{l}=0$.
But then
\[
   u_{d,a,l}u_{d,-a,l}=(1+z^da\ve{l})(1-z^da\ve{l})=1,
\]
completing the proof of~(a).
\\
(b) Let $u\in\Azs$ be a unit. Then $\lideal{u}=\Azs$ and thus $1\in A$ is a reduced
generator of $\lideal{u}$.
In~\cite[Cor.~4.13(a) and its proof]{GS02} it has been shown that the reduction of a
single polynomial in $\Azs$ can be described by left multiplication with suitable
elementary units.
In other words, there exist elementary units $u_1,\ldots,u_t\in\Azs$ such that
$1=u_t\cdot\ldots\cdot u_1u$ which proves the assertion.
\end{proof}

It should be noticed that from a coding theoretic point of view the elementary units
are not desirable if~$d$ is big.
Indeed, since the coefficients of $z,\,z^2,\ldots,z^{d-1}$ are zero, the same is
true for the coefficients of any component $u^{(l)}$ and thus the code
$\v(\lideal{u^{(l)}})$ has small distance.
This argument, of course, does not apply if $d=1$ and we will proceed with that more
specific case.
These units are not only candidates for the construction of good codes but, as we
will see next, will lead us to the existence of the desired minimal \scy\ codes.
To this end, we will now construct units whose $l$-th component have a prescribed
degree.

\begin{cor}\label{C-units}
Let $\sigma\in\AutF(A)$ and $l\in\{1,\ldots,r\}$ such that
$\sigma(\ve{l})\not=\ve{l}$.
Then we have
\begin{arabiclist}
\item For any $a\in A$ and any $i\in\N_0$ the element $u_a(i):=1+za\sigma^i(\ve{l})$
      is an elementary unit in $\Azs$. Its inverse is given by $u_{-a}(i)$.
\item For any $d\in\N_0$ and any units $a_1,\ldots,a_d$ in~$A$ the polynomial
      $u:=u_{a_1}(1)\cdot\ldots\cdot u_{a_d}(d)$ is a unit in $\Azs$ and
      satisfies $\deg_zu^{(l)}=d=\deg_zu$.
\end{arabiclist}
\end{cor}
\begin{proof}
(1) If $\deg_z u_a(i)=0$ the assertion is trivial.
Thus let us assume $\deg_z u_a(i)=1$.
Note that, with the notation of Lemma~\ref{L-eltunits}, $u_a(i)=u_{1,a,l'}$ where~$l'$ is such that
$\sigma^i(\ve{l})=\ve{l'}$.
From\eqnref{e-lorder} we know that $o_l(\sigma)=o_{l'}(\sigma)$
and by assumption this number is bigger than~$1$.
Thus $o_{l'}(\sigma)\nmid \deg_z u_a(i)$ and
Lemma~\ref{L-eltunits}(a) implies the assertion.
\\
(2) Without loss of generality let $d>0$.
Let $u:=u_{a_1}(1)\cdot\ldots\cdot u_{a_d}(d)$ where $a_1,\ldots,a_d$ are units
in~$A$.
From part~(a) we know that $u$ is a unit in $\Azs$ and has $\deg_z u\leq d$.
In order to show that $\deg_z u=d$ we compute the $z^d$-term of~$u$.
It is given by
\begin{align*}
   &\hspace*{-2em} \big(za_1\sigma(\ve{l})\big)\cdot\big(za_2\sigma^2(\ve{l})\big)\cdot\ldots\cdot
   \big(za_d\sigma^d(\ve{l})\big)\\[.7ex]
   =&\ z^d\big(\sigma^{d-1}(a_1)\sigma^{d-2}(a_2)\cdot\ldots\cdot\sigma(a_{d-1})a_d\big)
        \big(\sigma^d(\ve{l})\cdot\ldots\cdot\sigma^d(\ve{l})\big)\\[.7ex]
   =&\ z^da\sigma^d(\ve{l}),
\end{align*}
where $a:=\sigma^{d-1}(a_1)\sigma^{d-2}(a_2)\cdot\ldots\cdot\sigma(a_{d-1})a_d$.
Since $a_1,\ldots,a_d$ are units in~$A$ the same is true for~$a$.
Thus $a\sigma^d(\ve{l})\not=0$ and we have $\deg_z u=d$.
Finally, $\deg_z u^{(l)}=d$ since
$\ve{l}z^da\sigma^d(\ve{l})=z^da\sigma^d(\ve{l})\not=0$.
\end{proof}

We would like to mention that for the unit $u$ thus constructed $\deg_z u^{(l')}<d$
whenever $l'\not=l$.
This can easily be seen from the above.

The following theorem combines our findings about the existence of minimal
\scy\ \CC{}s.
The proof follows from Theorem~\ref{T-ideals}(4),
Lemma~\ref{L-compl0}, and Corollary~\ref{C-units}(2).

\begin{theo}\label{T-existminCCC}
Let $\sigma\in\AutF(A)$ and $l\in\{1,\ldots,r\}$.
Put $k:=\deg_x\pi_l$. Then the following are equivalent:
\begin{romanlist}
\item $\sigma(\ve{l})\not=\ve{l}$.
\item For any $d\in\N_0$ one can construct a minimal \scy\ \CC\ with parameters
      $(n,k,dk)$ and support~$\{l\}$.
      The Forney indices of the code are all equal to~$d$.
\item There exists a \scy\ \CC\ with nonzero complexity and support~$\{l\}$.
\end{romanlist}
\end{theo}

Notice that the considerations so far do not lead to any insight about the
quality of a minimal \scy\ \CC, that is, about the distance.
The following examples, however, suggest that this construction is worth being
investigated with respect to distance properties.
The codes given below are all optimal with respect to their distance.
As for the general situation, we wish to add that the codes constructed in
Theorem~\ref{T-existminCCC}(ii) are
{\em compact}, which in this case (rank~$k$ dividing the complexity) means
that the Forney indices are all the same~\cite[Cor.~4.3]{McE98}.
In general, compact codes are better candidates for good codes; for instance,
codes attaining the generalized Singleton bound\eqnref{e-MDS} are always
compact~\cite[Proof of Thm.~2.2]{RoSm99}.

\begin{exa}\label{E-minC3}
We begin with the case~$n=3$ over $\F:=\F_4=\{0,1,\alpha,\alpha^2\}$ where
$\alpha^2+\alpha+1=0$.
Thus $A=\F[x]/_{\DS\ideal{x^3-1}}$
and we have the prime factor decomposition $x^3-1=\pi_1\pi_2\pi_3$ where
$\pi_1=x+1,\;\pi_2=x+\alpha$, and $\pi_3=x+\alpha^2$.
The corresponding primitive idempotents are
\[
 \ve{1}=x^2+x+1,\ \ve{2}=\alpha x^2+\alpha^2x+1,\
 \ve{3}=\alpha^2x^2+\alpha x+1
\]
as can readily be seen by verifying
$(\ve{i}\;\mod\pi_j)=\delta_{ij}$ for $i,\,j=1,2,3$.
We will use the automorphism $\sigma\in\AutF(A)$ defined by
$\sigma(x)=x^2$.
One easily checks that $\sigma(\ve{2})=\ve{3}$ and vice versa.
Hence $\Pi_{\sigma}=(1)(2,3)$.
We will construct minimal \scy\ codes with support~$\{2\}$ by using the
construction of units in Corollary~\ref{C-units} for $l=2$.
Choose the units
\[
   v_1=u_1(1),\ v_2=u_{\alpha}(2),\ v_3=u_{\alpha^2}(3),\
   v_4=u_{\alpha}(4),\ v_5=u_{\alpha^2}(5),\ v_6=u_{\alpha}(6)\in\Azs
\]
and put $g^{(\delta)}:=\ve{2}(v_1\cdot\ldots\cdot v_{\delta})$ for $\delta=1,\ldots,6$.
From Corollary~\ref{C-units}(2) we know that $\deg_z g^{(\delta)}=\delta$ and that
$\cC^{(\delta)}:=\v(\lideal{g^{(\delta)}})$ is a \scy\ code with parameters
$(3,1,\delta)_4$.
We used a computer algebra program and computed the distances of these codes which
turn out to be very good in each case.
Indeed, the respective distances are
\[
  \dist(\cC^{(1)})\!=\!6,\,\dist(\cC^{(2)})\!=\!9,\,\dist(\cC^{(3)})\!=\!12,\,
  \dist(\cC^{(4)})\!=\!14,\,\dist(\cC^{(5)})\!=\!16,\,\dist(\cC^{(6)})\!=\!18.
\]
For $\delta=1,\ldots,5$ the distances attain the Griesmer bound\eqnref{e-Griesmer},
hence these codes are optimal (for $\delta=1,2,3$ this is even the generalized Singleton
bound\eqnref{e-MDS}).
For $\delta=6$ the computed distance is just one less than the Griesmer bound,
which in this case is~$19$.
It should be added that, as to our knowledge, it is unknown whether there exists
any code over~$\F_4$ with parameters $(3,1,6)$ and distance~$19$.
We think it is worth presenting these codes explicitly.
Recall from Theorem~\ref{T-ideals}(4) that $G^{(\delta)}:=\v(g^{(\delta)})$
is a generator matrix of~$\cC^{(\delta)}$.
These matrices are given by
\[
\begin{array}{ll}
   G^{(1)}=[z\!+\!1,\,\alpha z\!+\!\alpha^2,\,\alpha^2z\!+\!\alpha],&
   \!\!G^{(2)}=[\alpha z^2\!+\!z\!+\!1,\, z^2\!+\!\alpha z\!+\!\alpha^2,\,
            \alpha^2z^2\!+\!\alpha^2z\!+\!\alpha],\\[2ex]
  G^{(3)}=\begin{bmatrix}z^3\!+\!\alpha z^2\!+\!\alpha z\!+\!1\\
              \alpha z^3\!+\!z^2\!+\!\alpha^2 z\!+\!\alpha^2\\
               \alpha^2z^3\!+\!\alpha^2 z^2\!+\!z\!+\!\alpha\end{bmatrix}^{\!\!\sf T},&
  \!\!G^{(4)}=\begin{bmatrix}
            \alpha z^4\!+\!z^3\!+\!z^2\!+\!\alpha z\!+\!1\\
            z^4\!+\!\alpha z^3\!+\!\alpha^2z^2\!+\!\alpha^2 z\!+\!\alpha^2\\
            \alpha^2z^4\!+\!\alpha^2z^3\!+\!\alpha
            z^2\!+\!z\!+\!\alpha\end{bmatrix}^{\!\!\sf T},\\[4ex]
    G^{(5)}=\begin{bmatrix}
               z^5\!+\!\alpha z^4\!+\!\alpha z^3\!+\!z^2\!+\!z\!+\!1\\
               \alpha z^5\!+\!z^4\!+\!\alpha^2z^3\!+\!\alpha^2z^2\!+\!\alpha z\!+\!\alpha^2\\
               \alpha^2z^5\!+\!\alpha^2z^4\!+\!z^3\!+\!\alpha z^2\!+\!\alpha^2z\!+\!\alpha
            \end{bmatrix}^{\!\!\sf T}\!\!\!,&
  \!\!G^{(6)}=\begin{bmatrix}
               \alpha z^6\!+\!z^5\!+\!z^4\!+\!\alpha z^3\!+\!\alpha^2z^2\!+\!z\!+\!1\\
               z^6\!+\!\alpha z^5\!+\!\alpha^2z^4\!+\!\alpha^2z^3\!+\!\alpha z^2\!+\!\alpha z\!+\!\alpha^2\\
               \alpha^2z^6\!+\!\alpha^2z^5\!+\!\alpha z^4\!+\!z^3\!+\! z^2\!+\!\alpha^2z\!+\!\alpha
            \end{bmatrix}^{\!\!\sf T}\!\!\!.
\end{array}
\]
\end{exa}

\begin{exa}\label{E-minC5}
Now we consider the case $n=5$ over $\F=\F_4=\{0,1,\alpha,\alpha^2\}$.
In this case $x^5-1=\pi_1\pi_2\pi_3$ where $\pi_1=x+1,\,\pi_2=x^2+\alpha x+1$,
and $\pi_3=x^2+\alpha^2x+1$ and the corresponding
primitive idempotents are
\[
  \ve{1}=x^4\!+\!x^3\!+\!x^2\!+\!x\!+\!1,\;
  \ve{2}=\alpha x^4\!+\!\alpha^2 x^3\!+\!\alpha^2x^2\!+\!\alpha x,\;
  \ve{3}=\alpha^2x^4\!+\!\alpha x^3\!+\!\alpha x^2\!+\!\alpha^2x.
\]
We choose the automorphism defined via $\sigma(x)=x^2$.
Again it is easily seen that $\sigma(\ve{2})=\ve{3}$ and vice versa.
We will use Corollary~\ref{C-units} for $l=2$ in order to construct minimal
\scy\ codes with support~$\{2\}$.
We define
\[
   g^{(1)}:=\ve{2}u_1(1),\; g^{(2)}:=\ve{2}u_1(1)u_{\alpha}(2),\;
   g^{(3)}:=\ve{2}u_1(1)u_{\alpha}(2)u_{\alpha^2}(3).
\]
Then we know that $\deg_z g^{(m)}=m$ and that
$\cC^{(m)}:=\v(\lideal{g^{(m)}})$ is a \scy\ code over~$\F_4$ with parameters
$(5,2,2m)$ for $m=1,2,3$.
Again we computed the distances and they turn out to be optimal in each case.
In this case Theorem~\ref{T-ideals}(4) implies that the generator matrix of
$\cC^{(m)}$ is made up by the two rows
$\v(g^{(m)})$ and $\v(xg^{(m)})$.
They are computed as
\begin{align*}
   &G^{(1)}\!=\!\begin{bmatrix}
            0&\alpha+\alpha^2z&\alpha^2+\alpha z&\alpha^2+\alpha
            z&\alpha+\alpha^2z\\
            \alpha+\alpha z&\alpha^2 z&\alpha&\alpha^2+\alpha^2 z&\alpha^2+\alpha z
            \end{bmatrix},\\[2ex]
   &G^{(2)}\!=\!\begin{bmatrix}
            0\!&\!\alpha+\alpha^2z+\alpha^2z^2\!&\!\alpha^2+\alpha z+z^2\!&\!\alpha^2+\alpha
            z+z^2\!&\!\alpha+\alpha^2z+\alpha^2z^2\\
            \alpha+\alpha z+\alpha^2z^2\!&\!\alpha^2 z+z^2\!&\!\alpha+z^2\!&\!\alpha^2+\alpha^2 z+\alpha^2z^2\!&\!
            \alpha^2+\alpha z
            \end{bmatrix},\\[2ex]
   &{\scriptsize G^{(3)}\!=\!\begin{bmatrix}
            0&\alpha\!+\!z\!+\!\alpha^2z^2\!+\!\alpha^2z^3&\alpha^2\!+\!\alpha^2 z\!+\!z^2\!+\!\alpha z^3&
            \alpha^2\!+\!\alpha^2 z\!+\!z^2\!+\!\alpha z^3&\alpha\!+\!z\!+\!\alpha^2z^2\!+\!\alpha^2z^3\\
            \alpha\!+\!\alpha^2z\!+\!\alpha^2z^2\!+\!\alpha z^3&z\!+\!z^2\!+\!\alpha z^3&
            \alpha\!+\!z^2\!+\!\alpha^2z^3&\alpha^2\!+\!z\!+\!\alpha^2z^2&
            \alpha^2\!+\!\alpha^2z\!+\!\alpha^2z^3
            \end{bmatrix}.}
\end{align*}
The distances are $\dist(\cC^{(1)})=8,\ \dist(\cC^{(2)})=12$, and $\dist(\cC^{(3)})=16$,
which is in each case the Griesmer bound\eqnref{e-Griesmer} for codes over~$\F_4$ with
parameters $(5,2,2m)$.
\end{exa}

\begin{rem}\label{R-nonequiv}
In~\cite[Table~II]{GS03} some other sequences of codes over~$\F_4$ with
parameters $(3,1,\delta)$ for $\delta=1,\ldots,5$ and $(5,2,2m), m=1,2,3$ have been given.
They have the same distances as the ones given in the previous two examples, hence are
also optimal.
It is worth being pointed out that those codes and the ones presented here are {\em not\/}
strongly equivalent in the sense that we call
two codes $\im G$ and $\im G'$ {\em strongly equivalent\/} if
$G=G'P D$ where $P\in Gl_n(\F)$ is a permutation matrix and~$D\in Gl_n(\F)$ is a nonsingular
diagonal matrix.
In other words, codes are strongly equivalent if they differ only by a permutation
and a rescaling of the entries of the codewords.
Strongly equivalent codes have, of course, the same parameters and the same distance.
From a coding point of view they have the same properties and can
therefore be identified.
As a consequence, the two families of codes obtained in the examples above are significantly
different from those constructed earlier.
\end{rem}

\section{Orthogonal Sums of Minimal Cyclic Codes}
\setcounter{equation}{0}
In this section we will extend the existence result from Theorem~\ref{T-existminCCC}
to certain non  minimal \scy\ codes.
The main tool for this task is the orthogonality as stated in\eqnref{e-orthog}.
It leads directly to the following lemma.
This in turn will imply that the sum of minimal codes having pairwise
orthogonal generator polynomials is direct.
Again, let~$\F$ be a finite field such that~$|\F|$ and~$n$ are coprime and let
$\sigma\in\AutF(A)$ be a fixed automorphism.
We will make heavy use of the prime factor decomposition\eqnref{e-xn-1} and the
notations introduced in Definition~\ref{D-sigmapermut}.

\begin{lemma}\label{L-orthogunits}
Let $l_1,\ldots,l_t\in\{1,\ldots,r\}$ be such that $l_i\nequivs l_j$
for $i\not=j$.
Furthermore, put $I:=\{1,\ldots,r\}\backslash\{l\mid l\equivs l_i\text{ for some }i=1,\ldots,t\}$.
\begin{arabiclist}
\item Let $u\in\Azs$ be a unit with inverse $u^{-1}=\bar{u}$.
      Then
      \[
         \sum_{j\equivs l_i}u^{(j)}\sum_{j\equivs l_i}\bar{u}^{(j)}
         =\sum_{j\equivs l_i}\ve{j}
          \text{ for $i=1,\ldots,t$ and }
         \sum_{j\in I}u^{(j)}\sum_{j\in I}\bar{u}^{(j)}=\sum_{j\in I}\ve{j}.
      \]
\item For $i=1,\ldots,t$ let $u_i\in\Azs$ be a unit with inverse $u_i^{-1}=\bar{u}_i$
      and let $u\in\Azs$ be a unit with inverse $u^{-1}=\bar{u}$.
      Then the element $w:=\sum_{i=1}^t\sum_{j\equivs l_i}u_i^{(j)}+\sum_{j\in I}u^{(j)}$
      is a unit with inverse
      $w^{-1}=\sum_{i=1}^t\sum_{j\equivs l_i}\bar{u}_i^{(j)}+\sum_{j\in I}\bar{u}^{(j)}$.
\item Each polynomial $g\in\Azs$ with support
      $T_g=\{l_1,\ldots,l_t\}$ is reduced.
\end{arabiclist}
\end{lemma}
\begin{proof}
(1) The implication in\eqnref{e-orthog} yields
\[
  u\bar{u}=\sum_{i=1}^t\Big(\sum_{j\equivs l_i}u^{(j)}\sum_{j\equivs l_i}\bar{u}^{(j)}\Big)
   +\sum_{j\in I}u^{(j)}\sum_{j\in I}\bar{u}^{(j)}=1=\sum_{j=1}^r\ve{j}.
\]
From this the assertion follows immediately since the coefficients of each of the first~$t$
summands are contained in $\sum_{j\equivs l_i}\ve{j}A$ while those of the second sum are in
$\sum_{j\in I}\ve{j}A$ and all these sets are disjoint.
\\
(2) follows from~(1) along the same line of arguments.
\\
(3) Write $g=\sum_{i=1}^t g^{(l_i)}$.
By\eqnref{e-fkcoeff} the coefficients of~$z$ in  $g^{(l_i)}$ are contained
in $\sum_{l\equivs l_i}\ve{l}A$ for all $i=1,\ldots,t$.
But then no term of some component~$g^{(l_i)}$ can be right divisible by the
leading monomial of any other component.
\end{proof}

All this leads to the existence of units with prescribed degrees for
pairwise orthogonal components.

\begin{theo}\label{T-orthogminiCCC}
Let $l_1,\ldots,l_t\in\{1,\ldots,r\}$ be such that $l_i\nequivs l_j$ for $i\not=j$.
Furthermore assume $o_{l_i}(\sigma)>1$, that is,
$\sigma(\ve{l_i})\not=\ve{l_i}$, for all $i=1,\ldots,t$.
Then for all $d_1,\ldots,d_t\in\N_0$ there exists a unit $w\in\Azs$ such that
$g:=\sum_{i=1}^t w^{(l_i)}$ is reduced and $\deg_z w^{(l_i)}=d_i$ for
$i=1,\ldots,t$.
\end{theo}
\begin{proof}
From Corollary~\ref{C-units}(2) we know that for each $i=1,\ldots,t$ there exists a unit
$u_i$ such that $\deg_z u_i^{(l_i)}=d_i$.
Put $w:=\sum_{i=1}^t\sum_{j\equivs l_i}u_i^{(j)}+\sum_{i\in I}u_1^{(i)}$
where, again,
$I=\{1,\ldots,r\}\backslash\{l\mid l\equivs l_i\text{ for some }i=1,\ldots,t\}$.
Then Lemma~\ref{L-orthogunits}(2) and~(3) yield the desired results.
\end{proof}

Using Theorem~\ref{T-ideals}(4) we obtain immediately the existence of orthogonal sums
of minimal cyclic codes with prescribed Forney indices.

\begin{cor}\label{C-orthogSumCodes}
Let $l_1,\ldots,l_t\in\{1,\ldots,r\}$ be such that $l_i\nequivs l_j$ for $i\not=j$ and
such that $o_{l_i}(\sigma)>1$ for all $i=1,\ldots,t$.
Put $k_i:=\deg_x\pi_{l_i}$.
Then for all $d_1,\ldots,d_t\in\N_0$ there exists a \scy\ code $\cC\subseteq\F[z]^n$
with parameters $(n,k,\delta)$ where $k=\sum_{i=1}^t k_i$ and
$\delta=\sum_{i=1}^t k_id_i$. The support is given by $\{l_1,\ldots,l_t\}$.
\end{cor}

Note that, according to Theorem~\ref{T-ideals}(4), any \scy\ code with
support $\{l_1,\ldots,l_t\}$  has to have parameters of the type above.

The arguments above may be used to construct non-minimal codes with given parameters and
support consisting of indices with pairwise disjoint cycles directly out of minimal codes.
We formulate the result in terms of direct summands in $\F[z]^n$.

\begin{theo}\label{T-ConstrNonMini}
For $i=1,\ldots,t$ let $\cC_i\subseteq\F[z]^n$ be a minimal \scy\ code with support
$\{l_i\}$ and complexity~$\delta_i$ and assume $l_i\nequivs l_j$ for $i\not=j$.
Then $\cC:=\sum_{i=1}^t\cC_i\subseteq\F[z]^n$ is a \scy\ code, too.
Its rank is given by $\rank \cC=\sum_{i=1}^t\rank\cC_i=\sum_{i=1}^t\deg_x\pi_{l_i}$,
and its complexity is $\delta(\cC)=\delta_1+\ldots+\delta_t$.
Furthermore, $\cC=\oplus_{i=1}^t\cC_i$ and its Forney indices are given by
the union of the Forney indices of the codes $\cC_1,\ldots,\cC_t$.
\end{theo}
\begin{proof}
For all $i=1,\ldots,t$ let $\cC_i=\v(\lideal{g_i})$ where $g_i=u_i^{(l_i)}$ for
some unit $u_i\in\Azs$.
Put $g:=g_1+\ldots+g_t$ and $\cC:=\v(\lideal{g})$.
Then, by Lemma~\ref{L-orthogunits}(3), the polynomial~$g$ is reduced, and by
part~(2) of that lemma $g=\sum_{i=1}^t w^{(l_i)}$ for some suitable unit
$w\in\Azs$.
Hence, by Theorem~\ref{T-ideals}(3), the submodule~$\cC$ is a direct summand, and by
part~(4) of that theorem it it is the direct sum of $\cC_1,\ldots,\cC_t$ and has the
desired rank, complexity, and Forney indices.
\end{proof}

We wish to illustrate the above by an example indicating that this construction does
indeed lead to good codes.

\begin{exa}\label{E-length7overF8}
Let $n=7$ and $\F=\F_8=\{0,1,\alpha,\alpha^2,\ldots,\alpha^6\}$ where
$\alpha^3+\alpha+1=0$.
Then $x^7-1=\prod_{i=0}^6\pi_i,\ \text{ where } \pi_i=x-\alpha^i$.
Notice that, since all fields $K^{(i)}=\F[x]/_{\ideal{\pi_i}}$ are isomorphic to ~$\F_8$,
the automorphisms on $A=\F[x]/_{\ideal{x^7-1}}$ are fully determined by the
permutation $\Pi_{\sigma}$.
We choose the automorphism~$\sigma$ corresponding to the permutation
$\Pi_{\sigma}=(1,2)(3,4,5)(6)(7)$.
Moreover, we take the polynomials
$g_1=\ve{1}+z\ve{2}+z^2\ve{1}\alpha$ and $g_2=\ve{3}+z\ve{4}\alpha+z^2\ve{5}\alpha^2$.
Then $g_1=\ve{1}g_1$ and $g_2=\ve{3}g_2$.
Since both polynomials, being components, are reduced, Theorem~\ref{T-ideals}(4) tells us that
$\lideal{g_1}$ and $\lideal{g_2}$ are submodules of rank~$1$ and complexity~$2$ each.
It can be checked via some tedious but straightforward calculation that the associated matrices
$\v(g_i)$ are right invertible, thus both ideals are direct summands of $\Azs$.
Hence they are \scy\ codes over~$\F_8$ with parameters $(7,1,2)$ each.
Since $1\nequivs3$, the polynomial $g=g_1+g_2$ is reduced (see Lemma~\ref{L-orthogunits}(3))
and $\lideal{g}$ is a direct summand according to Theorem~\ref{T-ConstrNonMini}.
A minimal generator matrix of the code $\v(\lideal{g})\subseteq\F_8[z]^7$ is given by
\[
{\footnotesize
  \begin{bmatrix}\!
    1\!+\!z\!+\!\alpha z^2\!&\! 1\!+\!\alpha^6z\!+\!\alpha z^2\!&\! 1\!+\!\alpha^5z\!+\!\alpha z^2\!&\!
    1\!+\!\alpha^4z\!+\!\alpha z^2\!&\!1\!+\!\alpha^3z\!+\!\alpha z^2\!&\! 1\!+\!\alpha^2z\!+\!
    \alpha z^2\!&\! 1\!+\!\alpha z\!+\!\alpha z^2\\
    1\!+\!\alpha z\!+\!\alpha^2 z^2\!&\! \alpha^5\!+\!\alpha^5z\!+\!\alpha^5z^2\!&\! \alpha^3\!+\!
    \alpha^2z\!+\!\alpha z^2\!&\!\alpha\!+\!\alpha^6z\!+\!\alpha^4z^2\!&\! \alpha^6\!+\!
    \alpha^3z\!+\!z^2\!&\! \alpha^4\!+\!z\!+\!\alpha^3z^2\!&\!\alpha^2\!+\!\alpha^4z\!+\!\alpha^6z^2
  \!\end{bmatrix}\!\!.
}
\]
The first and second row generate the codes $\v(\lideal{g_1})$ and $\v(\lideal{g_2})$, respectively.
Again, all codes involved are optimal with respect to their distance.
Both the codes $\v(\lideal{g_i}),\,i=1,2,$ have distance~$21$,
which is the generalized Singleton bound\eqnref{e-MDS}.
Hence these codes are MDS codes in the sense of \cite{RoSm99}.
The code $\v(\lideal{g})$ has distance~$18$, which is the optimum value for codes over~$\F_8$ with
parameters $(7,2,4)$ due to the Griesmer bound\eqnref{e-Griesmer}.
\end{exa}

Finally we wish to comment on the existence of cyclic codes with arbitrary support.
We will briefly sketch that the existence result of Corollary~\ref{C-orthogSumCodes}
is not true without the assumption $l_i\nequivs l_j$ for $i\not=j$.
More precisely, in general it is not possible to arbitrarily prescribe the degrees
of the components of a reduced polynomial.
In order to see this, we consider a reduced polynomial~$g$
with support~$T_g$ containing at least two indices belonging to the same cycle
of~$\Pi_{\sigma}$.
Without restriction assume $S=\{1,\ldots,c\}\subseteq T_g$ and
$\sigma(\ve{i})=\ve{i+1}$ for all~$i=1,\ldots,o-1$ where
$o:=o_1(\sigma)\geq c$.
Let $\deg_zg^{(l)}=d_l$.
Then for $l=1,\ldots,c$ the highest coefficient of $g^{(l)}$ is in
$\sigma^{d_l}(\ve{l})A=\ve{(l+d_l-1\text{ mod } o) +1}A$
(the exponents arise from the fact that we have to compute modulo~$o$ with remainders
in $\{1,\ldots,o\}$ instead of $\{0,\ldots,o-1\}$).
Hence the reducedness of~$g$ implies that the numbers
\[
   (1+d_1),\ldots,(c+d_c)\text{ are pairwise different modulo $o$.}
\]
But for $c>1$ this puts a restriction on the degrees~$d_l$ of the components~$g^{(l)}$
(even without using the fact that~$g$ is the generator polynomial of
a code, i.~e., of a direct summand).
In case $c=o$, a second restriction arises if~$g$ generates a \scy\ code.
In that case not all~$d_l$ can be the same for otherwise one can easily see that~$g$
cannot be extended to a unit in~$\Azs$, see Corollary~\ref{C-dirSummands}.
It remains an open question whether there are further restrictions on the degrees
of the components.

\section{Open Problems}
We wish to close the paper with some problems open to future research.
As described at the end of the last section, in the general situation it remains open as to which
Forney indices (and complexity) a \scy\ code can attain.
But from a coding theoretic point of view an investigation of
\scy\ codes with respect to their distance is much more important.
More precisely, it needs to be investigated whether one can relate the distance of a \CCC\ to
some properties of the generator polynomial (or any other suitable generating polynomial
of the associated left ideal).
As a starting point one might begin with minimal codes.
In particular we think it is worth to investigate the construction of minimal codes
via units as described in Corollary~\ref{C-units}(2).
Furthermore, it is also unclear which automorphisms should be chosen for obtaining good codes.
Finally, the class of all cyclic codes of a given length needs to be investigated with
respect to strong equivalence in the sense given in Remark~\ref{R-nonequiv}.
First ideas can be found in~\cite{La03}, they indicate that one may restrict to certain
automorphisms in order to cover all equivalence classes.
A detailed positive result would considerably reduce the amount of data to be
investigated for the search of good cyclic codes.

\bibliographystyle{abbrv}
\bibliography{literatureAK,literatureLZ}

\begin{thebibliography}{10}

\bibitem{Fo70}
G.~D. Forney~Jr.
\newblock Convolutional codes {I}: {A}lgebraic structure.
\newblock {\em IEEE Trans. Inform. Theory}, 16:720--738, 1970.
\newblock (see also corrections in {\em IEEE Trans.\ Inf.\ Theory},
  vol.~17,1971, p.~360).

\bibitem{Fo75}
G.~D. Forney~Jr.
\newblock Minimal bases of rational vector spaces, with applications to
  multivariable linear systems.
\newblock {\em SIAM J. on Contr.}, 13:493--520, 1975.


\bibitem{GS02}
H.~Gluesing-Luerssen and W.~Schmale.
\newblock On cyclic convolutional codes.
\newblock Preprint 2002. Submitted. Available at http://front.math.ucdavis.edu/
  with ID-number RA/0211040.


\bibitem{GS03}
H.~Gluesing-Luerssen and W.~Schmale.
\newblock Distance bounds for convolutional codes and some optimal codes.
\newblock Preprint 2003. Submitted. Available at http://front.math.
  \mbox{ucdavis}.edu/ with ID-number RA/0305135.


\bibitem{Hu74}
T.~W. Hungerford.
\newblock {\em Algebra}.
\newblock Springer, New York, 1974.

\bibitem{JoZi99}
R.~Johannesson and K.~S. Zigangirov.
\newblock {\em Fundamentals of Convolutional Coding}.
\newblock IEEE Press, New York, 1999.

\bibitem{Ju73}
J.~Justesen.
\newblock New convolutional code constructions and a class of asymptotically
  good time-varying codes.
\newblock {\em IEEE Trans. Inform. Theory}, IT-19:220--225, 1973.

\bibitem{Ju75}
J.~Justesen.
\newblock Algebraic construction of rate $1/\nu$ convolutional codes.
\newblock {\em IEEE Trans. Inform. Theory}, IT-21:577--580, 1975.

\bibitem{La03}
B.~Langfeld.
\newblock Minimal cyclic convolutional codes.
\newblock Diploma Thesis at the University of Oldenburg (Germany). Available at
  http://www-m9.ma.tum.de/dm/homepages/ langfeld/thesis.pdf, 2003.

\bibitem{MS77}
F.~J. Mac{W}illiams and N.~J.~A. Sloane.
\newblock {\em The Theory of Error-Correcting Codes}.
\newblock North-Holland, 1977.

\bibitem{MCJ73}
J.~L. Massey, D.~J. Costello, and J.~Justesen.
\newblock Polynomial weights and code constructions.
\newblock {\em IEEE Trans. Inform. Theory}, IT-19:101--110, 1973.

\bibitem{McE98}
R.~J. Mc{E}liece.
\newblock The algebraic theory of convolutional codes.
\newblock In V.~Pless and W.~Huffman, editors, {\em Handbook of Coding Theory,
  Vol.~1}, pages 1065--1138. Elsevier, Amsterdam, 1998.

\bibitem{Pi75}
P.~Piret.
\newblock On a class of alternating cyclic convolutional codes.
\newblock {\em IEEE Trans. Inform. Theory}, 12:64--69, 1975.

\bibitem{Pi76}
P.~Piret.
\newblock Structure and constructions of cyclic convolutional codes.
\newblock {\em IEEE Trans. Inform. Theory}, 22:147--155, 1976.

\bibitem{Ro79}
C.~Roos.
\newblock On the structure of convolutional and cyclic convolutional codes.
\newblock {\em IEEE Trans. Inform. Theory}, 25:676--683, 1979.

\bibitem{Ro01}
J.~Rosenthal.
\newblock Connections between linear systems and convolutional codes.
\newblock In B.~Marcus and J.~Rosenthal, editors, {\em Codes, Systems, and
  Graphical Models}, pages 39--66. Springer, Berlin, 2001.

\bibitem{RSY96}
J.~Rosenthal, J.~M. Schumacher, and E.~V. York.
\newblock On behaviors and convolutional codes.
\newblock {\em IEEE Trans. Inform. Theory}, 42:1881--1891, 1996.

\bibitem{RoSm99}
J.~Rosenthal and R.~Smarandache.
\newblock Maximum distance separable convolutional codes.
\newblock {\em Appl.\ Algebra Engrg.\ Comm.\ Comput.}, 10:15--32, 1999.

\bibitem{SGR01}
R.~Smarandache, H.~Gluesing-Luerssen, and J.~Rosenthal.
\newblock Constructions of {MDS}-convolutional codes.
\newblock {\em IEEE Trans. Inform. Theory}, 47(5):2045--2049, 2001.

\bibitem{Ve85}
M.~Ventou.
\newblock Automorphisms and isometries of some modular algebras.
\newblock In {\em Algebraic algorithms and error-correcting codes; {P}roc. 3rd
  {I}nternational {C}onf.~{AAECC}-3}, pages 202--210. Springer Lecture Notes in
  Computer Science~229, 1985.

\end{thebibliography}
\end{document}